\documentclass[a4paper]{article}
\usepackage{amssymb}

\newcommand{\bc}{\begin{center}}
\newcommand{\ec}{\end{center}}
\newcommand{\be}{\begin{equation}}
\newcommand{\ee}{\end{equation}}
\newcommand{\bea}{\begin{eqnarray}}
\newcommand{\eea}{\end{eqnarray}}
\newcommand{\ba}{\begin{array}}
\newcommand{\ea}{\end{array}}

\newcommand{\ra}{\rightarrow}
\def\i{{\bf 1}\!\!{\rm I}}
\def\p{\psi}
\def\l{\lambda}

\def\t{\tau}
\def\s{\sigma}
\def\a{\alpha}
\def\G{\Gamma}
\def\n{\nu}
\def\f{\varphi}
\def\o{\otimes}
\def\x{\infty}
\def\R{\mathbb{R}}
\def\N{\mathbb{N}}
\def\A{\mathbb{A}}
\def\V{\mathbf{V}}
\def\W{\mathbf{W}}

\def\e{\varepsilon}
\begin{document}

\begin{center}
{\large {\bf ON INFINITE DIMENSIONAL QUADRATIC VOLTERRA
OPERATORS}}
\footnote{The work supported by NATO-TUBITAK PC-B programme}\\[1.1cm]

{FARRUKH MUKHAMEDOV}\footnote{E-mail: E-Mail: far75m@yandex.ru; Corresponding author}\\
Department of Mechanics and Mathematics,\\
National University of Uzbekistan,\\
Vuzgorodok, 700174,\\
Tashkent, Uzbekistan, \\[2mm]
 {HASAN AKIN}\footnote{E-mail: akinhasan@harran.edu.tr} \ \ \
{SEYIT TEMIR} \footnote{E-mail:
temirseyit@harran.edu.tr}\\[2mm]

Department of Mathematics, \\
Arts and Science Faculty, \\Harran University, 63200, \\
\c{S}anliurfa, Turkey
\end{center}

\begin{abstract}
In this paper we study a class of quadratic operators named by
Volterra operators on infinite dimensional space. We prove that
such operators have infinitely many fixed points and the set of
Volterra operators forms a convex compact set. In addition, it is
described its extreme points. Besides, we study certain limit
behaviors of such operators and give some more examples of
Volterra operators for which their trajectories do not converge.
Finally, we define a compatible sequence of finite dimensional
Volterra operators and prove that any power of this sequence
converges in weak topology.

\end{abstract}

\section{Introduction}

It is known that the theory of Markov processes is a
well-developed field of mathematics which has various applications
in physics, biology and so on. But there are some physical models
which cannot be described by such processes. One of such models is
a model related to population genetics. Namely, consider a
biological population, that is,a community of organisms closed
with respect to reproduction \cite{B}. Assume that every
individual in this population belongs to precisely one of the
species $1,2,\cdots,n$. The scale of species is such that the
species of parents $i$ and $j$ unambiguously determine the
probability of every species $k$ for the first generation of
direct descendants. We denote this probability (the heredity
coefficient) by $p_{ij,k}$. It is obvious that $p_{ij,k}\geq 0$
and $\sum\limits_{k=1}^np_{ij,k}=1$ for all $i,j$. Assume that the
population is so large that frequency fluctuations can be
neglected. Then the state of the population can be described by
the tuple $x=(x_1,x_2,\cdots,x_n)$ of species probabilities, that
is, $x_i$ is the fraction of the species $i$ in the population.In
this case of panmixia (random interbreeding), the parent pairs $i$
and $j$ arise for a fixed state $x=(x_1,x_2,\cdots,x_n)$ with
probability $x_ix_j$. Hence \be \label{int}
x'_k=\sum_{i,j=1}^np_{ij,k}x_ix_j \ee is the total probability of
the species $k$ in the first generation of direct descendants. The
set $ S^{n-1}=\{x\in \R^n :\  x_i\geq 0, \
\sum\limits_{i=1}^{n}x_i=1\}$ is an $(n-1)$-dimensional simplex.
Since, $x_k'\geq 0$ and $\sum\limits_{i=1}^{n}x'_i=1$, the {\it
quadratic stochastic operator} defined by formula (\ref{int}) maps
$S^{n-1}$ into itself. In this setting an evolution of the system
is described by this  operator acting on the  simplex. Note that
the notion of quadratic operator firstly introduced by Bernstein
in \cite{B}. To investigations of such kind operators devoted a
lot papers (see \cite{L2} for review). One of the central problem
in this theory is to study limit behavior of quadratic  operators
(see \cite{U}).

In \cite{K},\cite{L1},\cite{SG},\cite{V} the authors investigated
limit behavior and ergodic properties of trajectories of the
quadratic stochastic operators. But these operators do not occupy
quantum systems, so it is natural to investigate quantum quadratic
operators. In \cite{GM1},\cite{GM2} a notion of quantum quadratic
stochastic operators defined on von Neumann algebra has been
introduced. It includes as a particular case of quadratic
stochastic operators. In \cite{GM2},\cite{M2} some ergodic and
stability properties of such operators were studied. But it would
be more interesting to investigate one of the simplest case in
which that operators  act on infinite dimensional algebras.

In this paper we are going to consider quadratic operators on
infinite dimensional commutative algebra. In this setting an
infinite dimensional simplex is not weak compact, therefore in
general we cannot state that every quadratic operator has at least
one fixed point. This is an infinite dimensionality phenomenon. We
will  study a class of quadratic operators named by Volterra
operators. The paper is organized as follows. In section 2 we give
some preliminary on quadratic operators defined on von Neumann
algebra and describe a form of such operators defined on
$\ell^\x$. Besides, we demonstrate an example of quadratic
operator which has no fixed points. In section 3 we define
quadratic Volterra operators and study its certain properties. In
particular, we show that such operators have infinitely many fixed
points. In section 4 we prove that  the set of Volterra operators
forms a convex compact set and describe its extreme points. In the
next section 5 we study certain limit behaviors of such operators
and give some more examples of Volterra operators for which their
trajectories do not converge. Finally, in the last section 6 we
define a compatible sequence of finite dimensional Volterra
operators and prove that any power of this sequence converges in
weak topology. It should be noted that finite dimensional Volterra
operators were studied in \cite{G}.

Note that a part of the results have been announced in \cite{M3}.

\section{Preliminary and quadratic operators}

Let us recall some definitions. Let $B(H)$ be the algebra of
linear bounded operators on a separable Hilbert space $H$.  Let
$M\subset B(H)$ be a von Neumann algebra with unit $\i$. By $M_+$
we denote the set of all positive elements of $M$. Weak (operator)
closure of algebraic tensor product ${M}\odot{M}$ in $B(H\otimes
H)$ is denoted by ${M}\otimes{M}$, and it is called {\sl tensor
product} of  $M$ into itself. For detail we refer a reader to
\cite{BR}.

By $S(M)$ and $S(M\otimes M)$ it is denoted the set of all normal
states on $M$ and $M\otimes M$ respectively. Let $U:M\o M\ra M\o
M$ be a linear operator such that $U(x\o y)=y\o x$ for all $x,y\in
M$.

{\bf Definition 2.1.}\cite{M1} A linear operator $P:M\to M\o M$ is
said to be  {\it quantum quadratic stochastic operator (q.q.s.o.)}
if it is normal and  satisfies the following conditions:
\begin{enumerate}
   \item[(i)] $P\i_{M}=\i_{M\o M}$, where $\i_{M}$ and $\i_{M\o M}$ are units of
   algebras $M$ and ${M\o M}$ respectively;
   \item[(ii)] $P(M_{+})\subset (M\o M)_{+}$;
   \item[(iii)] $UPx=Px$ for every  $x\in {M}$.\\
\end{enumerate}

Define an operator $\tilde V:S(M\otimes M)\to S(M)$ as follows
\be\label{1} \tilde V(\tilde\f)(x)=\tilde \f(Px), \ \ \tilde \f\in
S(M\otimes M), \ x\in M. \ee

The operator $\tilde V$ is called {\it conjugate quadratic
operator (c.q.o.)}. Further for the shortness instead of $\tilde
V(\f\otimes\p)$ we will  write  $\tilde V(\f,\p)$, where $\f,\p\in
S(M)$. Note that the relation (iii) implies that \be\label{sym}
\tilde V(\f,\p)=\tilde V(\p,\f). \ee

In \cite{M2} we have proved that every c.q.o. uniquely defines
q.q.s.o. Therefore it is enough to consider c.q.o.

By means of $\tilde V$ one can define an operator $V:S(M)\to S(M)$
by \be\label{qo} V(\f)=\tilde V(\f,\f), \ \ \f\in S(M),\ee
which is called {\it quadratic operator (q.o.)}.\\

{\bf Observation 2.1.} Here we give how linear operator and
q.q.s.o. related with each other. Let $T:M\to M$ be a linear
positive normal operator (i.e. $Tx\geq 0$ whenever $x\geq 0$) such
that $T\i=\i$. Define a linear operator $P:M\to M\o M$ as follows
\be\label{qqso} Px=\frac{Tx\o\i+\i\o Tx}{2}, \ \ \ x\in M. \ee It
is clear that $P$ is q.q.s.o. Then associated c.q.o. and q.o. have
the following form respectively: \bea\label{qqso1} \tilde
V(\f,\p)(x)&= &\frac{1}{2}(\f+\p)(Tx), \nonumber \\
V(\f)(x) &= & \f(Tx), \ \ \ x\in M,\eea for every $\f,\p\in S(M)$.
Thus linear operator can be viewed as a particular case of
q.q.s.o. If $T$ is the identity operator, then from (\ref{qqso1})
we can find that the associated q.o. also would be the identity
operator of $S$. The set of all q.q.s.o. associated with linear
operators we denote
by ${\cal QL}(M)$. \\

In the paper we are going to consider a case when the von Neumann
algebra $M$ is a infinite-dimensional commutative discrete
algebra, i.e.
$$
M=\ell^{\x}=\{x=(x_{n}) : \ x_n\in \R, \
\|x\|_{\x}=\sup\limits_{n\in \N}|x_i| \},
$$
then the set of all normal functionals defined on $\ell^\x$
coincides with
$$
\ell^{1}=\{x=\{x_n\} :\  \|x\|_1=\sum_{k=1}^{\x}|x_k|<\x \}
$$
(i.e. $\ell^1$ is a pre-dual space to $\ell^\x$, namely
$(\ell^1)^*=\ell^\x$) and $S(\ell^{\x})$ with
$$
S=\{x=(x_n)\in \ell^1 : x_i\geq 0, \sum_{n=1}^{\x}x_n=1 \}.
$$
It is known \cite{R} that  $S=\overline{convh(Extr S)}$, where
$Extr(S)$ is the extremal points of $S$ and $convh(A)$ is the
convex hall of a set $A$.

Any extremal point $\f$ of $S$ has the following form
$$
\f=(\underbrace{0,0,...,1}_n,0,...),
$$
for some $n\in\N$.  Such elements will be denoted  as $e^{(n)}$.\\

The following Theorem describes c.q.o. when $M=\ell^\x$.

 {\bf Theorem 2.1.} {\it Every c.q.o.
$\tilde V$ defines an infinite dimensional matrix
$(p_{ij,k})_{i,j,k\in \N}$ such that \be\label{qso} p_{ij,k}\geq
0, \ \ p_{ij,k}=p_{ji,k}, \ \ \sum_{k=1}^{\infty}p_{ij,k}=1, \ \
i,j\in N. \ee Conversely, every such matrix defines c.q.o. $\tilde
V$ as follows: \be\label{qso1} (\tilde
V(x,y))_k=\sum_{i,j=1}^{\infty}p_{ij,k}x_iy_j, \ \ k\in \N, \
x=(x_i), y=(y_i)\in S. \ee }

{\bf Proof.} Let $\tilde V$ be a c.q.o. For every
$e^{(n)},e^{(m)}\in Extr(S)$ put
$$
p_{mn,k}=(\tilde V(e^{(m)},e^{(n)}))_k, \ \ m,n,k\in \N.
$$
According to positivity of $e^{(n)},n\in \N$ and (ii) (see
def.1.1) we get $p_{mn,k}\geq 0$. It follows from (\ref{sym}) that
$ \tilde V(e^{(m)},e^{(n)})=\tilde V(e^{(n)},e^{(m)}),$ which
implies that $p_{mn,k}=p_{nm,k}$. Since $\tilde
V(e^{(m)},e^{(n)})\in S$ we find
$\sum\limits_{k=1}^{\infty}p_{mn,k}=1.$ Note that we have
$$
(\tilde V(x,y))_k=\sum_{i,j=1}^{\infty}p_{ij,k}x_iy_j, \ \ k\in
\N.
$$
for every $x=(x_i),y=(y_i)\in S$.

Conversely, let $(p_{ij,k})$ be a matrix satisfying (\ref{qso}).
Define $P:\ell^\x\to\ell^\x\o\ell^\x$ as follows
$$
(Pf)_{ij}=\sum_{k=1}^{\x}p_{ij,k}f_k, \ \ \ \ i,j\in\N,
$$
for every $f=(f_k)\in\ell^\x$. The condition (\ref{qso}) implies
that $P$ is a q.q.s.o. In particulary, we have \be\label{qso2}
Pe^{(k)}=\sum_{i,j\in\N}p_{ij,k}e^{(i)}\o e^{(j)}. \ee Let $\tilde
V$ be the c.q.o. associated with $P$. Take arbitrary $x,y\in S$.
Then using (\ref{qso2}) we find
$$
(\tilde V(x,y))_k=x\o y(Pe^{(k)})=\sum_{k=1}^{\x}p_{ij,k}x_iy_j,
$$
here $x\o y =(x_iy_j)\in S(\ell^\x\o\ell^\x)$.

Thus the theorem is proved.

We note that in this case q.o. $V$ defined by (\ref{qo}) has the
following form: \be\label{qo2} (V(x))_k=
\sum_{i,j=1}^{\infty}p_{ij,k}x_ix_j \ \ k\in \N, \ x=(x_i)\in S.
\ee

The constructed matrix $(p_{ij,k})_{i,j,k\in\N}$ is called {\it
determining matrix of q.o. $V$}.

{\bf Observation 2.2.} Let $T:\ell^{\x}\to\ell^{\x}$ be a positive
identity preserving operator. Then it is easy to see that this
operator can be represented as infinite dimensional stochastic
matrix $(p_{ij})_{i,j\in\N}$, i.e. $ p_{ij}\geq 0$,
$\sum\limits_{j=1}^{\x}p_{ij}=1$ for every $i,j\in\N$.

Then the determining matrix $(p_{ij,k})_{i,j,k\in\N}$
corresponding to q.o.  given by (\ref{qqso1}) is defined as
$$
p_{ij,k}=\frac{p_{ik}+p_{jk}}{2}, \ \ \ i,j,k\in\N.\\
$$

{\bf Observation 2.3.} It is known that the set $S$ is not compact
in norm topology of $\ell^1$, even in
$\s(\ell^1,\ell^\infty)$-topology. This is the difference between
finite and infinite dimensional cases. In  finite dimensional case
 every q.o. $V:S^{n-1}\to S^{n-1}$ has at least one fixed point
(i.e. $V(x)=x$, $x\in S^{n-1}$). In the infinite dimensional
setting, not every q.o. has fixed points. Indeed, define a linear
operator $T:\ell^\infty\to \ell^\infty$ as follows
$$
T(x_1,x_2,\cdots,x_n,\cdots)=(x_2,\cdots,x_{n+1},\cdots),
$$
$(x_n)\in \ell^\infty$. It is clear that $T$ is positive and
$T\i=\i$. Now consider q.q.s.o. defined by (\ref{qqso}). Then by
Observation 2.1 q.o. $V$ acts as follows
$$
V(\f_1,\f_2,\cdots,\f_n,\cdots)=(0,\f_1,\f_2,\cdots,\f_n,\cdots)
$$
where $(\f_n)\in S$. It is easy to see that this operator has no
fixed points belonging to $S$. \\

\section{Volterra operators}

In this section we define Volterra operators and give some their
properties.

Recall that a convex set $C\subset S $ is called {\sl face}, if
$\l x+(1-\l)y\in C$, where $x,y\in S$ ш $\l\in (0,1)$, implies
that $x,y\in C$. For $\f,\p\in S$ denote
$\G(\f,\p)=\{\l\f+(1-\l)\p : \l\in [0,1] \}$.

{\bf Definition 3.1.} An operator $V$ defined by (\ref{qo}) is
called {\it Vollterra operator} if $\tilde V(\f,\p)\in \G(\f,\p)$
is valid for every $\f,\p\in Extr(S)$.

By ${\cal QV}$ we denote the set of all quadratic operators
defined on $S$, and the set of all Volterra operators is denoted
by ${\cal V}$.

{\bf Proposition 3.1.} {\it Let  $V\in {\cal QV}$ be a q.o. Then
$V$ is Volterra if and only if  the determining matrix
$(p_{ij,k})$ of this operator satisfies the following property:
\be\label{her} p_{ij,k}=0, \ \ \textrm{if} \ \ k\notin\{i,j\}.
\ee}

{\bf Proof.} Let $V$ be a Volterra operator. Then from definition
3.1 we infer that
$$
\tilde V(e^{(i)},e^{(j)})= p_{ij,i}e^{(i)}+p_{ij,j}e^{(j)}.
$$
This yields that $p_{ij,i}+p_{ij,j}=1$, so (\ref{her}) is valid.
The converse implication easily follows from Theorem 2.1. The
proposition is proved.

Note that the condition (\ref{her}) biologically means that each
individual can inherit only the species of the parents.

From Theorem 2.1 and Proposition 3.1 we immediately get the
following

{\bf Proposition 3.2.} {\it Let  $V_1,V_2\in {\cal V}$ be two
Volterra operators such that for every $e^{(i)},e^{(j)}, \ i,j\in
\N$ the equality holds $ \tilde V_1(e^{(i)},e^{(j)})= \tilde
V_2(e^{(i)},e^{(j)})$, then $V_1=V_2$.}

{\bf Theorem 3.3.} {\it Let $V\in {\cal QV}$ be a q.o. Then $V$ is
Volterra operator if and only if it can be represented as follows:
\be\label{vol} (V(x))_{k}=x_k(1+\sum_{i=1}^{\x}a_{ki}x_i), \ \
k\in \N, \ee where \be\label{anti}a_{ki}=-a_{ik}, \ |a_{ki}|\leq 1
\ \ \textrm{for every} \ \ k,i\in\N\ee.}

{\bf Proof.} From Definition 3.1 and Proposition 3.1.  one gets
$p_{kk,k}=1, k\in \N$. Then from (\ref{qo2}) we obtain \bea
(V(x))_{k}&=&\sum_{i,j=1}^{\x}p_{ij,k}x_ix_j\nonumber\\
&= & p_{kk,k}x^{2}_{k}+\sum_{i=1, i\neq
k}p_{ik,k}x_ix_k+\sum_{j=1, j\neq k}p_{kj,k}x_kx_j, \ \ k\in\N,
\nonumber \eea whence keeping in mind  $p_{ij,k}=p_{ji,k}$ we
infer that
$$
(V(x))_{k}=x_k(1+2\sum_{i=1, i\neq k}^{\x}p_{ik,k}x_i), \ \ k\in
\N.
$$
Using  $\sum\limits_{i=1}^{\x}x_i=1$ we have
$$
(V(x))_{k}=x_k(1+\sum_{i=1, i\neq k}^{\x}(2p_{ik,k}-1)x_i), \ \
k\in\N.
$$
Setting  $a_{ki}=2p_{ik,k}-1$ while $i\neq k$, and  $a_{kk}=0$, it
yields (\ref{vol}). The inequality  $0\leq p_{ik,k}\leq 1$ implies
that $|a_{ki}|\leq 1$. Taking into account $p_{ik,k}+p_{ik.i}=1,$
we have
$$
a_{ki}+a_{ik}=2p_{ik,k}-1+2p_{ki,i}-1=2(p_{ik,k}+p_{ik.i}-1)=0.
$$
Therefore  $a_{ki}=-a_{ik}$.

The converse implication is obvious. This completes the proof.

{\bf Corollary 3.4.} {\it Let $V\in {\cal QV}$ be a q.o. Then $V$
is Volterra operator if and only if $\tilde V$ can be represented
as follows: \be\label{vol1} (\tilde
V(x,y))_{k}=\frac{1}{2}\left(x_k(1+\sum_{i=1}^{\x}a_{ki}y_i)+
y_k(1+\sum_{i=1}^{\x}a_{ki}x_i)\right), \ \ k\in \N. \ee }

Recall that an element $x\in S$ is called fixed point of $V$ if
$V(x)=x$. The set of all fixed points of $V$ is denoted by
$Fix(V)$. For given subset $K$ of $\N$ set
$$
S^K=\{x\in S \ :\ x_i=0, \ \forall i\in\N\setminus K \}.
$$

{\bf Corollary 3.5.} {\it For every Volterra operator $V$ the
following assertions hold:}
\begin{enumerate}
   \item[(i)] {\it every face of $S$ invariant with respect to
   $V$;}
   \item[(ii)] $ Extr(S)\subset Fix(V)$.
\end{enumerate}

The proof immediately follows from Theorem 3.3 since every face of
$S$ is $S^K$ for some $K\subset\N$ and $\{e^{(i)}\}=S^{\{i\}}$ for
every $e^{(i)}\in Extr(S)$.

Put
$$
ri S^K=\{x\in S^K | x_i>0, \ \ \forall i\in K \}.
$$

{\bf Corollary 3.6.} {\it Let $V$ be a Volterra operator, then the
relation  holds $ V(ri S^K)\subset riS^K$ for every $K\subset\N$}

{\bf Proof.}  Let $x_k>0, k\in K$, then according to the equality
$a_{kk}=0$ and (\ref{vol}) we heve \bea
(V(x))_{k}&=&x_k(1+a_{k1}x_1+...+a_{k,k-1}x_{k-1}+a_{k,k+1}x_{k+1}+...) \nonumber\\
&\geq & x_k(1-x_1-...-x_{k-1}-x_{k+1}-...)=x^{2}_{k}>0.  \nonumber
\eea The corollary is proved.

{\bf Remark 3.1.} From Theorem 3.3 we see that the identity
operator $Id:S\to S$, i.e.
$$
(Id(x))_k=x_k, \ \ \ k\in\N
$$
is Volterra operator. From  Proposition 3.2, Observations 2.1 and
2.2 we infer that ${\cal QL}(l^{\x})\cap{\cal V}=Id$.

{\bf Theorem 3.7.} {\it Let $V\in {\cal V}$ be a Volterra
operator, then it is a bijection of $S$.}

{\bf Proof.} Let us first show that $V$ is injective. Assume that
there are two elements $x,y\in S (x\neq y)$ such that
\be\label{eq} V(x)=V(y)\ee

Without loss of generality we may assume that $x_i>0, y_i>0,
\forall i\in \N.$ If it is not true then there is a face $S^K$,
for some subset $K\subset\N$, of $S$ such that $x,y\in ri S^K$,
i.e. $x_i>0,y_i>0$, $\forall i\in K$. According to Corollaries 3.5
and 3.6 we have $V(S^K)\subset S^K$, therefore we may restrict
$V$ to $S^K$. From  (\ref{eq}) one gets that
$$
x_k(1+\sum_{i=1}^{\infty}a_{ki}x_i)=y_k(1+\sum_{i=1}^{\infty}a_{ki}y_i),
$$
or \be\label{eq1}
(x_k-y_k)(1+\sum_{i=1}^{\infty}a_{ki}y_i)=-x_k\sum_{i=1}^{\infty}a_{ki}(x_i-y_i).
\ee

We have
$$
1+\sum_{i=1}^{\infty}a_{ki}y_i\geq
1-y_1-y_2-...-y_{k-1}-y_{k+1}-...=y_k>0,
$$
whence $x_k>0$ with (\ref{eq1}) implies  that \be\label{eq2} sgn
(x_k-y_k)=-sgn \sum_{i=1}^{\infty}a_{ki}(x_i-y_i). \ee Hence
$$
(x_k-y_k)\sum_{i=1}^{\infty}a_{ki}(x_i-y_i)\leq 0, \ \ k\in \N,
$$
whence
$$
\sum_{k=1}^{\infty}(x_k-y_k)\sum_{i=1}^{\infty}a_{ki}(x_i-y_i)\leq
0.
$$
Note that the last series absolutely converges, since  \bea
|\sum_{k=1}^{\infty}(x_k-y_k)\sum_{i=1}^{\infty}a_{ki}(x_i-y_i)|&\leq
&
\sum_{k=1}^{\infty}|x_k-y_k|\sum_{i=1}^{\infty}|a_{ki}||x_i-y_i|\nonumber \\
&\leq &
\sum_{k=1}^{\infty}(x_k+y_k)\sum_{i=1}^{\infty}(x_i+y_i)=4<\infty.
\nonumber \eea According to $a_{ki}=-a_{ik}$ we find
$$
\sum_{k=1}^{\infty}(x_k-y_k)\sum_{i=1}^{\infty}a_{ki}(x_i-y_i)= 0.
$$
Consequently,
$$
(x_k-y_k)\sum_{i=1}^{\infty}a_{ki}(x_i-y_i)= 0, \ \ k\in \N.
$$
The equality (\ref{eq2}) with the last equality imply that $x=y$.
Thus,  $V:S\to S$ is injective.

Now let us show that $V$ is onto. Denote
$$
A_1=\{ [1,n]\subset \N : n\in \N \}, \ \  A_2=\{ a\subset [1,n] :
|[1,n]\setminus a|\geq 2,  n\in \N \},
$$
$$
A_3=\{ b\subset \N : a\subset b, \ a\in A_1\cup A_2, \
|\N\setminus b|<\infty,  \},
$$
$$ A=A_1\cup A_2\cup A_3.$$
Order $A$ by inclusion, i.e.  $a\leq b$ means that $a\subset b$
for $a,b\in A$. It is clear that $A$ is a completely ordering set.
To prove that $V$ is surjective we will use transfer induction
method with respect to the set $A$. Obviously, that for the first
element $\{1\}$ of the set  $A$, the operator $V$ on $S^{\{1\}}$
is surjective (see Corollary 3.5 and \cite{G}). Assume that for an
element $a\in A$ the operator $V$ is surjective  on $S^b$ for
every $b<a$. Let us show that it is surjective on $S^a$. Suppose
that $V(S^a)\neq S^a$. For the boundary  $\partial S^a$ of $S^a$
we have $\partial S^a=\bigcup\limits_{c\in A: c<a}S^c$. According
to the assumption of the induction one gets \be\label{eq3}
V(\partial S^{a})=\partial S^{a}. \ee On the other hand,  there
exist $x,y\in ri S^a$ such that $x\in V(S^a), \ y\notin V(S^a)$.
The segment $[x,y]$ contains at least one boundary point $z$ of
the set $V(S^a)$. Since $V:S^a\to V(S^a)$ is continuous and
bijection, then the boundary point goes to boundary one. Therefore
for $z\in ri S^a$, $V^{-1}(z)\in
\partial S^a$, which contradicts to (\ref{eq3}). Thus the theorem
is proved.

\section{The set of Volterra operators}

In this section we will prove that the set ${\cal V}$ is compact.

Now endow ${\cal QV}$ with a topology which is defined by  the
following system of semi-norms:
$$
p_{\f,\p,k}(V)=|(V(\f,\p))_k|, \ \ \  V\in {\cal QV},
$$
where $\f,\p\in S$ and $k\in\N$. This topology is called {\it weak
topology} and is denoted by $\t_w$.

A net $\{V_{\n}\}$ of quadratic operators converges to $V$ with
respect to the defined topology if for every $\f,\p\in S$ and
$k\in\N$
$$
(V_{\n}(\f,\p))_k\to (V(\f,\p))_k
$$
is valid.

Since ${\cal V}\subset {\cal QV}$, therefore on  ${\cal V}$ we
consider the induced topology  by ${\cal QV}$.

We note that in \cite{M1} we have  proved that the set of all
quantum quadratic stochastic operators defined on semi-finite von
Neumann algebra, without normality condition, forms is a weak
compact convex set. In the present situation we cannot apply the
mentioned result since our q.q.s.o. are normal. In general, the
set of all normal q.q.s.o. is not weak compact\footnote{Each state
$\omega\in S(M)$ defines a linear positive operator as
$T(x)=\omega(x)\i$. So according to Observation 2.1 the set of all
normal states can be included to ${\cal QV}$. Therefore, we can
consider the induced weak topology(defined as above) on $S(M)$. It
is clear that this topology coincides with $*$-topology on $S(M)$,
but in this topology $S(M)$ is not compact. Hence, ${\cal QV}$ is
not weak compact.}. Therefore, here we use another method to prove
that ${\cal V}$ is weak compact.

Denote the set of all matrices $(a_{ki})$ satisfying (\ref{anti})
by $\A$. It is clear that $\A$ is convex. The set $\A$ can be
considered as a subset of the space
$$
\ell^{\x}(\N\times\N)=\{x=(x_{n,m}) : x_{n,m}\in \R, \ n,m\in\N, \
\|x\|_{\x}=\sup\limits_{n,m\in \N}|x_{n,m}| \}.
$$

It is well-known \cite{BR} that the space
$$
\ell^{1}(\N\times\N)=\{x=(x_{n,m}) : x_{n,m}\in \R, \ n,m\in\N, \
\|x\|_1=\sum_{n,m\in\N}|x_{n,m}|< \infty\}
$$
is pre-dual to $\ell^{\x}(\N\times\N)$, i.e.
$\ell^{1}(\N\times\N)^*=\ell^{\x}(\N\times\N)$. Therefore on
$\ell^{\x}(\N\times\N)$ we can consider
$\sigma(\ell^{\x}(\N\times\N),\ell^{1}(\N\times\N))$-topology. In
the sequel we will denote it as $\t$.  According to Alaoglu-Banach
theorem the set $\A$ is
$\sigma(\ell^{\x}(\N\times\N),\ell^{1}(\N\times\N))$-weak compact
in $\ell^{\x}(\N\times\N)$. From  Theorem 3.3 we conclude that
every $(a_{ki})$ matrix with the property (\ref{anti}) defines  a
Volterra operator $V$ of the form (\ref{vol}) ( see also
(\ref{vol1})). So, it is defined a map $T:\A\to{\cal V}$. It is
clear that Theorem 3.3 and Proposition 3.2 imply that this map is
bijection and convex.

{\bf Theorem 4.1.} {\it The map $T:(\A,\t)\to ({\cal V},\t_w)$ is
continuous.}

{\bf Proof.} Let a net $(a^{(\n)}_{ki})\subset \A$ converge to
$(a_{ki})$ in the weak topology. This means that for an arbitrary
$\e>0$ and  every $k,i\in\N$ there is $\n_0(ki)$ such that
$|a^{(\n)}_{ki}-a_{ki}|<\e$ for every $n\geq\n_0(ki)$. Denote
$V^{(\n)}=T((a^{\n}_{ki}))$ and $V=T((a_{ki}))$.

Take any $x,y\in S$. Then there is a number $N_0\in\N$ such that
\be\label{es} \sum_{i=N_0+1}^{\x}x_i<\e, \ \ \ \ \
\sum_{i=N_0+1}^{\x}y_i<\e. \ee

Now consider two separate cases.

{\tt Case (i)}. In this case we assume that $1\leq k\leq N_0$.
Then according to Corollary 3.4 and using (\ref{anti}),(\ref{es})
we infer that \bea |(\tilde V^{(\n)}(x,y))_k-(\tilde V(x,y))_k|
&\leq & \frac{1}{2}\left(\sum_{i\in
\N}(y_kx_i+x_ky_i)|a^{(\n)}_{ki}-a_{ki}|\right) \nonumber \\
&\leq
&\frac{1}{2}\left(\sum_{i=1}^{N_0}(x_i+y_i)|a^{(\n)}_{ki}-a_{ki}|\right)+
\sum_{i=N_0+1}^{\x}(x_i+y_i)< 3\e\nonumber \eea for every
$\n\geq\max\{\n_0(ki):k,i\leq N_0\}$. Here we have used that
$\sum\limits_{i=1}^{N_0}(x_i+y_i)\leq
\sum\limits_{i=1}^\x(x_i+y_i)=2$.

{\tt Case (ii)}. Now assume that $k\geq N_0+1$. Using the above
argument we have \bea |(\tilde V^{(\n)}(x,y))_k-(\tilde
V(x,y))_k|&\leq & \frac{1}{2}\left(\sum_{i\in
\N}(y_kx_i+x_ky_i)|a^{(\n)}_{ki}-a_{ki}|\right)\nonumber \\
&\leq & y_k+x_k\leq \sum_{i=N_0+1}^{\x}(x_i+y_i)< 2\e\nonumber
\eea for every $\n\geq\max\{\n_0(ki):k,i\leq N_0\}$. Thus the map
$T$ is continuous. The theorem is proved.

{\bf Corollary  4.2.} {\it The set ${\cal V}$ is weak convex
compact.}

The proof immediately comes from that $\A$ is compact and $T$ is
continuous.

We say that q.o. $V\in{\cal QV}$ is {\it pure} if  for every
$\f,\p\in S$ the relation holds
$$
\tilde V(\f,\p)\in Extr \G(\f,\p) =\{\f,\p\}.
$$

It is clear that pure q.o. are Volterra.

{\bf Proposition 4.3.} {\it The set ${\cal V}$ is convex.
Moreover, $V$ is extreme  point of ${\cal V}$ if and only if it is
pure.}

{\bf Proof.} Convexity of ${\cal V}$ is obvious. Let $V$ be a pure
q.o. Let us assume that there exits $\l\in (0,1)$ and operators
$V_1,V_2\in {\cal V}$ such that $V=\l V_1+(1-\l)V_2$.

Let $\f\,\p\in Extr(S)$, then we have \be\label{conv} \tilde
V(\f,\p)=\l\tilde V_1(\f,\p)+(1-\l)\tilde V_2(\f,\p). \ee Without
loss of generality we may suppose that $\tilde V(\f,\p)=\f$, since
$V$ is pure. Therefore, the extremity of $\f$ with (\ref{conv})
implies that $V_i(\f,\p)=\f, \ i=1,2$. Hence,
$V_1(\f,\p)=V_2(\f,\p)$ for every $\f,\p\in Extr(S)$. According to
Proposition 3.2 one gets $V=V_1=V_2$. Thus $V\in Extr({\cal V})$.

Now let $V\in Extr({\cal V})$. Show that $V$ is pure. Assume that
$V$ is not pure, i.e. there is $\f_0,\p_0\in Extr(S)$ and a number
$\l\in (0,1)$ such that $ \tilde V(\f_0,\p_0)=\l\f_0+(1-\l)\p_0.$
Define q.o. $V_1$ and  $V_2$ as follows:
 \bea
\tilde V_1(\f_0,\p_0)& = &\f_0, \ \ \tilde V_2(\f_0,\p_0)=\p_0 \nonumber \\
 V_i(\f,\p) &=&  V(\f,\p) \ \ \forall \f,\p\in Extr(S),\f,\p\notin\{\f_0,\p_0\} .
\nonumber \eea Then again using Proposition 3.2 we get $ V=\l
V_1+(1-\l)V_2,$ which contradicts to the extremity of $V$. This
completes the proof.

We note that Proposition 4.3 can be also proved by means of
Theorem 3.3 and Corollary 3.4.

From Corollary 4.2 and Proposition 4.3 we have the following

{\bf Corollary  4.4.} {\it A Volterra operator $V\in{\cal V}$ is
extremal if and only if for the associated skew-symmetric matrix
$(a_{ki})$ the equality holds $|a_{ki}|=1$, for every $ k,i\in\N.$
}

The proof comes from that the extremal points of $\A$ satisfy the
last condition and the map $T$ is convex and bijection.

\section{A limit behavior of Volterra operators}

In this section we give some  limit theorems concerning
trajectories of Volterra operators.

 Let  $V:S\to S$ be a Volterra operator. Then according to
Theorem 3.3 it  has the form (\ref{vol}).

Denote \be\label{Q} \mathbf{Q}=\{y\in S : \
\sum_{i=1}^{\x}a_{ki}y_i\leq 0, \ \ k\in \N\}.\ee

It is clear that $\mathbf{Q}$ is convex subset of $S$.

{\bf Proposition 5.1.} {\it  For every Volterra operator $V$ the
relation holds $ \mathbf{Q}\subset Fix(V).$}

{\bf Proof.} Let $y\in \mathbf{Q}$ then \be\label{eq4}
(V(y))_{k}=y_{k}(1+\sum_{i=1}^{\x}a_{ki}y_i)\leq y_k, \ \ k\in
\N.\ee According to the equality
$\sum\limits_{i=1}^{\x}y_i=\sum\limits_{i=1}^{\x}(V(y))_i=1$, from
(\ref{eq4}) we find $(V(y))_{k}=y_k$ for every $k\in \N$, i.e.
$Vy=y$.

{\bf Theorem 5.2.}{\it Let $V$ be a Volterra operator such that
$\mathbf{Q}\neq\emptyset$. Suppose  $x^0\in riS$ (i.e.
$x^{0}_{i}>0, \forall i\in \N$) such that $Vx^0\neq x^0$ and the
limit  $\lim\limits_{n\to\x}V^nx^0$ exits. Then
$\lim\limits_{n\to\x}V^nx^0\in \mathbf{Q}$.}

{\bf Proof.} Let  $x^0\in riS$ and
$\lim\limits_{n\to\x}x^{(n)}=\tilde x,$ where $x^{(n)}=V^nx^0, \
n\in \N$. Denote $\tilde x=(q_1,q_2,...,q_n,...)$. It is clear
that  $V\tilde x=\tilde x$. Hence \be\label{eq5}
q_k=q_k(1+\sum_{i=1}^{\x}a_{ki}q_i), \ \ k\in N. \ee Set
$I_+=\{i\in \N | q_i>0\}, \ I_0=\{i\in N | q_i=0\}$. If $k\in
I_+$, then from (\ref{eq5}) we get
$$
\sum_{i=1}^{\x}a_{ki}q_i=0, \ \ k\in I_+.
$$
Assume that there is  $k_0\in I_0$ such that
$$
\sum_{i=1}^{\x}a_{k_{0}i}q_i>0.
$$
Since $x^{(m)}_k\to q_k$, then there is $m_0\in \N$ such that
\be\label{eq6} \sum_{i=1}^{\x}a_{k_{0}i}x^{(m)}_{i}>0, \ \
\mbox{for every} \ \ m\geq m_0. \ee

According to Corollary 3.6 we have $x^{(m)}\in riS, \ \forall m\in
\N$, i.e. $x^{(m)}_{k}>0, \ \forall m,k\in \N$. The inequality
(\ref{eq6}) with one
$$
x^{(m+1)}_{k_{0}}=x^{(m)}_{k_{0}}(1+\sum_{i=1}^{\x}a_{k_{0}i}
x^{(m)}_{i})>x^{(m)}_{k_{0}}  \ \ \forall m\geq m_0
$$
implies  that $x^{(m+1)}_{k_{0}}>x^{(m)}_{k_{0}}$, which
contradicts to $x^{(m)}_{k_{0}}\to q_{k_{0}}=0$. Therefore if
$k\in I_0$, then $\sum\limits_{i=1}^{\x}a_{ki}q_i\leq 0$. Thus
$\tilde x\in \mathbf{Q}$. The theorem is proved.

Given $V$ Volterra operator and  $K\subset\N$. Set $V_K=V|_{S^K}$.
Let  $\mathbf{Q}_K$ be the set $\mathbf{Q}$ corresponding to
$V_K$. Then from Theorem 5.2 and Corollary 3.6 we immediately get

{\bf Corollary 5.3.} {\it Let $\mathbf{Q}_K\neq\emptyset$ and
$x^0\in riS^K$ (i.e. $x^{0}_{i}>0, \forall i\in K$) such that
$Vx^0\neq x^0$ and the limit $\lim\limits_{n\to\x}V^nx^0$ exists.
Then $\lim\limits_{n\to\x}V^nx^0\in \mathbf{Q}_K$.}

{\bf Corollary 5.4.} {\it If a Volterra operator $V$ has an
isolated fixed point $x^0\in Fix(V)$(i.e. there is a weak neighbor
$U(x^0)\subset S$ of $x^0$ such that $U(x^0)\cap Fix(V)=\{x^0\}$)
such that $x^0\in riS$. Then for any $x\in riS$,$x\notin Fix(V)$
the limit $\lim\limits_{n\to\x}V^nx$ does not exists. }

{\bf Proof.} Assume that $\lim\limits_{n\to\x}V^nx=\bar x$ exists.
Then according to Theorem 5.2 we have $\bar x\in \mathbf{Q}$.
Since $x^0\in Fix(V)$, $x^0\in riS$ imply that $x^0\in\mathbf{Q}$.
Convexity of $\mathbf{Q}$ yields that $\l\bar
x+(1-\l)x^0\in\mathbf{Q}$ for every $\l\in[0,1]$. But this
contradicts the fact that $x^0$ is isolated. This completes the
proof.

{\bf Remark 5.1.} It is known \cite{G} that the set $\mathbf{Q}$
is not empty for any Volterra operator in finite dimensional
setting. But unfortunately, in our situation $\mathbf{Q}$ can be
empty.

Let us give some more examples of q.o. for which $\mathbf{Q}$ is
empty and non empty.

{\bf Example 5.1.} Let us consider a Volterra operator defined as
follows: \bea \left\{ \ba{ll}
(Vx)_{2k-1}=x_{2k-1}(1-a^{(k)}x_{2k}), \\[2mm]
(Vx)_{2k}=x_{2k}(1+a^{(k)}x_{2k-1}), \ k\in \N \ea \right.
\nonumber \eea уфх $a^{(k)}>0$ ш $|a^{(k)}|\leq 1$.

Let us describe $\mathbf{Q}$ for defined $V$. To this end we
should find solutions of the system:
$$
\left\{ \ba{ll}
-a^{(k)}x_{2k}\leq 0, \\[2mm]
a^{(k)}x_{2k-1}\leq 0, \ k\in \N \ea \right.
$$

One easily gets that  $\mathbf{Q}=\{x\in S : x_{2k-1}=0, \ k\in\N
\}$. So $\mathbf{Q}\neq\emptyset$.

Let $x\in riS$, then the trajectory of $x$ is defined as the
following recurrent relations
 \bea \left\{ \ba{ll}
x_{2k-1}^{(m+1)}=x_{2k-1}^{(m)}(1-a^{(k)}x_{2k}^{(m)}), \\[2mm]
x_{2k}^{(m+1)}=x_{2k}^{(m)}(1+a^{(k)}x_{2k-1}^{(m)}), \\
\ea \right. k\in\N, m\in N.\nonumber \eea

According to $a^{(k)}>0$ we find $1+a^{(k)}x_{2k-1}^{(m)}>0$,
hence we have $ x^{(m+1)}_{2k}\geq x^{(m)}_{2k},$ therefore
$\{x^{(m)}_{2k}\}$ is non-decreasing sequence. From $0\leq
1-a^{(k)}x_{2k}^{(m)}\leq 1$ it follows that $\{x^{(m)}_{2k-1}\}$
is non-increasing sequence, such that $ 0\leq
x^{(m)}_{2k},x^{(m)}_{2k-1}\leq 1.$ So the limits
$$
\lim_{m\to\x}x^{(m)}_{2k-1}=\a_{2k-1}, \ \
 \lim_{m\to\x}x^{(m)}_{2k}=\beta_{2k}
$$
exits.

According to Theorem 5.2 we infer that $\alpha_{2k-1}=0$ for every
$k\in\N$.

Now let $x\notin riS$, then denote $I_x=\{ k\in\N : x_{k}=0\}$.
Then using Corollary 3.5 we find $V(S^{\N\setminus
{I_x}})=S^{\N\setminus I_x}$. The restriction of  $V$ to
$S^{\N\setminus I_x}$ is denoted by  $V_{\N\setminus I_x}$. From
definition of $S^{\N\setminus I_x}$ we find that $x\in ri
S^{\N\setminus I_x}$, whence according to Corollary 5.3 we obtain
$$
\lim_{m\to\x}x^{(m)}_{2k-1}=0, \ \ \lim_{m\to\x}x^{(m)}_{2k}=
\left\{ \ba{ll}
\beta_{2k}, \ \ 2k\in \N\setminus I_x,\\[2mm]
0, \ \ 2k\in I_x.\\
\ea \right.
$$

{\bf Example 5.2.}  Let us define a Volterra operator as follows:
\be\label{exm2} (V(x))_{k}=x_k(1+\sum_{i=1}^{\x}a_{ki}x_i), \ \
k\in \N, \ee where $a_{ki}=(-1)^i$, $a_{ik}=-a_{ki}$ at $i\geq
k+1$.

Then it is not hard to check that the set  $\mathbf{Q}$ consists
of the solutions the following system \bea\label{exm21} \left\{
\begin{array}{lllll}
\sum\limits_{k=2}^{\x}(-1)^{k+1}x_{k}\leq 0, \\
x_{1}+\sum\limits_{k=3}^{\x}(-1)^{k}x_{k}\leq 0,  \\
-x_{1}+x_{2}+\sum\limits_{k=4}^{\x}(-1)^{k+1}x_{k}\leq 0,  \\
.......................................................... \\
\sum\limits_{k=2}^{n-1}(-1)^{n+k}x_{k}+\sum_{k=n+1}^{\x}(-1)^{n+k+1}x_{k}\leq
0. \\
..........................................................
\end{array}
\right. \eea Whence one gets  $x_{n}\leq x_{n+1}$ for every
$k\in\N$. Since  $x_{1}\geq 0$ and $x_n\to 0$ at $n\to\x$, we
obtain $x_n=0, \ \forall n\in\N$, which is impossible, because of
$\sum\limits_{k=1}^{\x}x_k=1$. Consequently,
$\mathbf{Q}=\emptyset$.

Now let us look for the set $Fix(V)$.  Let  $x^0\in riS$, i.e.
$x^{0}_{k}>0, \ \forall k\in \N$, be a fixed point of $V$. It
follows from (\ref{exm2}),(\ref{exm21}) that
$$
x^{0}_{1}=x^{0}_{2}=...x^{0}_{k}=..., \ \ k\in \N,
$$
but this equality is impossible since $x^{0}_{1}\neq 0$ and
$x^0_n\to 0$. Hence, inner fixed points for $V$ does not exit. So
there is a subset $I\subset\N$ such that $ I=\{k\in \N |
x^{0}_{k}=0\}.$ The set $\N\setminus I$ is finite. Indeed, assume
that $|\N\setminus I|=\x$, then consider a face $S^{\N\setminus
I}$. Then according to Corollary 3.5 $V_{\N\setminus I}$ is a
Volterra operator. It is clear that a point $ x^{0,\N\setminus
I}=\{x^{0}_{k} | \ k\in \N\setminus I\}$ is a fixed point of
$V_{\N\setminus I}$. From  (\ref{exm2}) and using the same
argument as above we find that the set $J=\{k\in \N |
x^{0,\N\setminus I}_{k}=0\}$ is non-empty, which contradicts to
the choice of $I$. Consequently, we infer that all fixed points of
$V$ lie on the faces $S^I$ such that $|\N\setminus I|<\x$. Thus we
conclude that the set $\mathbf{Q}$ turns out to be empty while the
set $Fix(V)$ is not. Therefore, Theorem 5.2 implies that if  $x\in
riS$ then the limit  $\lim\limits_{n\to\x}V^nx$ does not exists.
Now let $x\notin riS$, then for the set $I_x$ there are two
possibilities. The first case. Let $|\N\setminus I_x|=\x$, then
$x\in riS^{\N\setminus I_x}$. From condition $|\N\setminus
I_x|=\x$ analogously reasoning as above one can show that the set
$\mathbf{Q}_I$ is empty. According to Corollary 5.3 we infer that
the limit  $\lim\limits_{n\to\x}V^nx$ does not exist. The second
case. In this setting  $|\N\setminus I_x|<\x$, then the operator
$V_I$ reduces to finite dimensional operator, therefore the set
$\mathbf{Q}_I$  is not empty (see \cite{G}). So the limit
$\lim\limits_{n\to\x}V^nx$ exists since $|a_{ik}|=1$ (see
\cite{G}).\\

Now we will give a sufficient condition for $V$ which ensures that
the set $\mathbf{Q}$ is not empty.

Let  $V:S\to S$ be a Volterra operator which has the form
(\ref{vol}). Let $A=(a_{ki})$ be the corresponding skew-symmetric
matrix. Further we will assume that $A$ acts on $\ell^1$. A matrix
$A$ is called {\it finite dimensional} if $A(\ell^1)$ is finite
dimensional. We say that $A$ is  {\it finitely generated} if there
are a sequence of finite dimensional matrices $\{A_n\}$ such that
$\sup_n\|A_n\|<\infty$ and
$$
A=A_1\oplus A_2\oplus\cdots\oplus A_n\oplus\cdots
$$

\indent{\bf Proposition 5.5.} {\it Let  $A=(a_{ki})$ be the
skew-symmetric matrix corresponding  to a Volterra operator (see
(\ref{vol}), is finitely generated. Then the system
\be\label{ineq} \sum_{i=1}^{\x}a_{ki}y_i\geq 0, \ \ k\in \N \ee
has at least one element belonging to $S$.}

{\bf Proof.} First assume that $A$ is finite - dimensional, i.e.
there is  $n\in N$ such that $A(\ell^1)=\R^n$. According to
skew-symmetricity of  $A$ we find that $a_{ij}=0$ at  $i,j\geq
n+1$. Therefore we may assume that $A$ acts on $\R^n$. Then
(\ref{ineq}) is rewritten as follows \be\label{ineq1}
\sum_{j=1}^{n}a_{kj}y_j\geq 0, \ \ \ k=1,...,n. \ee According to
\cite{G} this system has a solution $y=\{y_k\}_{k=1}^{n}\in
S^{n-1}$ such that (\ref{ineq1}) holds. Now define an element
$\tilde y=\{\tilde y_k\}_{k=1}^{\infty}\in S$ as follows
$$
\tilde y_k=
\left\{
\ba{ll}
y_k, \ \ \textrm{if} \ \ 1\leq k\leq n\\
0, \ \ \textrm{if} \ \  k\geq n+1\\
\ea \right.
$$
It is evident that  $A\tilde y\geq 0$.

Now let us assume that $A$ is finitely generated, i.e.
$A=A_1\oplus A_2\oplus\cdots\oplus A_n\oplus\cdots$. Since
operators $A_n$ are finite dimensional, therefore suppose that for
every $n\in N$ there is $m_n\in \N$ such that  $A_n$ acts on
$\R^{m_n}$, i.e. $A_n:\R^{m_n}\to \R^{m_n}$. Consider the system
$$
A_ny^{(n)}\geq 0, \ \ \ n\in \N.
$$
According to the above argument, for every  $n\in \N$,  there is
an element  $z^{(n)}\in S^{m_n-1}$ such that $A_nz^{(n)}\geq 0$.
Define $z=(z_k)_{k=1}^{\infty}$ by
$$
z=\frac{1}{2}z^{(1)}\oplus\frac{1}{2^2}z^{(2)}\oplus...\oplus\frac{1}{2^n}z^{(n)}\oplus... \ .
$$
From \bea
\sum_{k=1}^{\infty}z_k=\frac{1}{2}\sum_{k=1}^{m_1}z^{(1)}_{k}+
\frac{1}{2^2}\sum_{k=1}^{m_2}z^{(2)}_{k}+...+\frac{1}{2^n}\sum_{k=1}^{m_n}z^{(n)}_{k}+...= \nonumber\\
=\frac{1}{2}+\frac{1}{2^2}+...+\frac{1}{2^n}+...=1.\nonumber \eea
we see that $z\in S$. The element $z$ is a solution of
(\ref{ineq}). Since
$$
Az=\frac{1}{2}A_1z^{(1)}\oplus...\oplus\frac{1}{2^n}A_nz^{(n)}\oplus...\geq 0.
$$
The proposition is proved.\\

{\bf Corollary 5.6.} {\it Let the condition of the previous
proposition is valid. Then the set $\mathbf{Q}$ is not empty.}

The proof immediately comes from Proposition 5.5 by changing the
matrix $A$ to $-A$, since  $-A$ is also skew-symmetric.

\section{Extension of finite dimensional Volterra operators}

In this section we are going to construct infinite dimensional
Volterra operators by means of finite dimensional ones.

Let $K_n=[1,n]\cap\N$ for every $n\in\N$. Consider a sequence
$V_{n]}:S^{K_n}\to S^{K_n}$ of finite dimensional Volterra
operators, i.e.
\be\label{seq}(V_{n]}(x))_{k}=x_{k}(1+\sum^{n}_{i=1}a_{ki}^{n]}x_{i})
\ \ k=1,\cdots,n, \ \ n\in\N, \ee here $(a_{ki}^{n]})$ is a
skew-symmetric matrix.

We say that this sequence of Volterra operators is {\it
compatible} if
\be\label{comp1}V_{n+1]}\upharpoonright{S^{K_n}}=V_{n]}\ee for
every $n\in \N$. The compatibility condition with  (\ref{seq})
implies that \be\label{comp2} a_{ki}^{n+1]}=a_{ki}^{n]}, \ \ \
\forall k,i\in\{1,\cdots,n\}.\ee

Denote
$$
S^{[n}=\{x=(x_{n},x_{n+1},\cdots): \ x_k\geq 0, \forall k\geq n,
\sum_{k=n}^\infty x_k=1\}, \ \ n\in\N.
$$

Let $\{W_{[n}:S^{[n}\to S^{[n}: n\in\N\}$ be a sequence of
Volterra operators
\be\label{seq1}(W_{[n}(x))_{k}=x_{k}(1+\sum^{\infty}_{i=n}a_{ki}^{[n}x_{i})
\ \ k\geq n, \ \ n\in\N. \ee

Define a sequence $\{\W_n:S\to S,\ n\in\N\}$ of  infinite
dimensional operators as follows
 \bea\label{vol4}
(\W_{n}(x))_{k} = \left\{
\begin{array}{ll}
(V_{n]}(x))_{k}, \ \  \textrm{if} \ \  n \leq k,\\[2mm]
 (W_{[n+1}(x))_{k}, \ \ \ \ \  \textrm{if} \ \ k \geq n+1,
\end{array}
\right.  \ \ \ n\in\N. \eea

According to Theorem 3.3 the defined operators are  Volterra.

{\bf Theorem 6.1.} {\it The sequence $\{\W_n\}$ of Volterra
operators weakly converges to a Volterra operator $\W$. Moreover,
if $V_{n]}$ are pure then $\W$ is so.}

{\bf Proof.} Let $x\in S$. If there is a finite subset $K$ of $\N$
such that $x\in S^K$, then according to the compatibility
condition (\ref{comp1}) we get  $\W(x)=\W_n(x)$ for all $n\geq
\max\{m: m\in K\}$.

Now assume that $x_i>0$ for all $i\in\N$. Let us prove that
$\{\W_{n}(x)\}$ is a Cauchy sequence with respect to weak
topology. Let $\e>0$ be an arbitrary number. Since $x\in S$ there
is  a number $n_0\in\N$ such that \be\label{n0}
\sum_{j=n+1}^{\infty}x_j<\e, \ \ \ \forall n\geq n_0 \ee

Consider several cases:

{\tt Case (i).} Suppose that $1\leq k\leq n$. Using (\ref{vol4}),
(\ref{comp2}),(\ref{seq}),(\ref{anti}) and (\ref{n0})  we have
\bea\label{cauchy} |(\W_{n}(x))_{k}-(\W_{n+p}(x))_{k}|& =
&|(V_{n]}(x))_{k}-(V_{n+p]}(x))_{k}|\nonumber
\\
& = &\left|x_{k}\bigg(1+\sum^{n}_{i=1}a_{ki}^{n]}x_{i}\bigg)-
x_{k}\bigg(1+\sum^{n+p}_{j=1}a_{kj}^{n+p]}x_{j}\bigg)\right|\nonumber
\\
&\leq&
x_{k}\left(\sum^{n+p}_{j=n+1}x_{j}\right)\leq\sum_{j=n+1}^{\infty}x_j<\e
\eea for all $n\geq n_0$.

{\tt Case (ii).} Assume $n+1 \leq k\leq n+p$. It then follows from
(\ref{seq1}),(\ref{vol4}) that \bea\label{cauchy2}
|(\W_{n}(x))_{k}-(\W_{n+p}(x))_{k}| & = &
|(W_{[n+1}(x))_{k}-(V_{n+p]}(x))_{k}|\nonumber \\
& = &
\left|x_{k}\bigg(1+\sum^{\infty}_{j=n+1}a_{kj}^{[n}x_{j}\bigg)-
x_{k}\bigg(1+\sum^{n}_{j=1}a_{kj}^{n+p]}x_{j}\bigg)\right|\nonumber\\
&\leq & x_{k}\sum^{\infty}_{j=1}|\gamma_{ki}|x_{j} \leq 2x_{k}<2\e
\eea for all $n\geq n_0$. Here \bea \gamma_{kj}= \left\{
\begin{array}{lll}
a_{kj}^{n+p]} \ \ \textrm{if} \  \  j\leq n,\\[2mm]
a_{kj}^{n+p]}+a_{kj}^{[n} \ \ \textrm{if} \ \ n+1\leq j\leq
n+p\\[2mm]
a_{kj}^{[n+1} \ \ \textrm{if}\  \  j\geq n+p+1\\[2mm]
\end{array}
\right.\nonumber  \eea

{\tt Case (iii).} Now assume that $k\geq n+p+1$ then from
(\ref{vol4}) we have \bea\label{cauchy4}
|(W_{[n}(x))_{k}-(W_{[n+p}(x))_{k}|& = &
\left|x_{k}\bigg(1+\sum^{\infty}_{j=n+1}a_{ki}^{[n+1}x_{j}\bigg)
-x_{k}\bigg(1+\sum^{\infty}_{j=n+p+1}a_{kj}^{[n+p+1}x_{j}\bigg)\right|\nonumber
\\
& \leq & 2x_{k}\sum^\infty_{j=n+1}x_{j}< 2\e^2 \eea

Hence the sequence  $(\W_{n}(x))$ is  Cauchy, therefore
$\W_{n}(x)\rightarrow \W(x)$.  By the same way we can show that
$\tilde{\W}_{n}(x,y)\rightarrow \tilde\W(x,y)$. Because of
$\W_{n}(e^{(i)},\ e^{(j)})\in \Gamma(e^{(i)},\ e^{(j)})$ and the
compatibility condition we find that $\W$ is Volterra. According
to (\ref{comp1}), for every $e^{(i)}$ and $e^{(j)}$ there is
$n_0\in\N$ such that $\tilde\W(e^{(i)},e^{(j)})=\tilde
V_{n]}(e^{(i)},e^{(j)})$ for all $n\geq n_0$.  Now if $V_{n]}$ is
pure for all $n\in\N$ then $\V$ is also pure. The theorem is
proved.

Let $\{\V_n:S\to S,\ n\in\N\}$ be a sequence of operators
associated with (see (\ref{vol4}))
$$
(W_{[n}(x))_k=x_k, \ \ \ k\geq n,\ n\in\N.
$$
According to Theorem 6.1 the defined sequence $\{\V_n\}$ converges
to a Volterra operator $\V$.

Now naturally comes a question: are the operators $\V$ and $\W$
equal? Next theorem gives an affirmative answer to this question.

{\bf Theorem 6.2.} {\it The operators $\V$ and $\W$ are equal.}

{\bf Proof.} Let $\e>0$ be an arbitrary number and $x\in S$ be
fixed. To prove the assertion it is enough to show for every
$k\in\N$ the relation holds
$$
|(\W_n(x))_k-(\V_n(x))_k|<\e.
$$

There is  a number $n_0\in\N$ such that (\ref{n0}) holds. Consider
two cases.

{\tt case (i).} Let $1\leq k\leq n$. Then  (\ref{vol4}) implies
that
$$
|(\W_n(x))_k-(\V_n(x))_k|=0,
$$
for every $n\geq n_0$.

{\tt case (ii).} Let $k\geq n+1$, then it follows from (\ref{n0})
that \bea |(\W_n(x))_k-(\V_n(x))_k|& =
&\left|x_{k}\bigg(1+\sum^{\infty}_{j=n+1}a_{kj}^{[n}x_{j}\bigg)-
x_{k}\right|\leq  x_k\sum^\infty_{j=n+1}x_j<\e.\nonumber
 \eea

Hence, we have proved the desired relation. This completes the
proof.

Thus according to the last Theorem we will consider only the
sequence $\{\V_n\}$. Now we are interested about the convergence
of powers of the sequence $\{\V_n\}$.

Let $V$ be an arbitrary Volterra operator. By $V^m$ we will denote
m-th iteration of $V$, i.e.
$V^m(x)=\underbrace{V(V\cdots(V}_m(x))\cdots).$ Before going to
formulate the result we need the following

{\bf Lemma 6.3.} {\it Let $V$ be an arbitrary Volterra operator.
Then $ (V^m(x))_k\leq 2^mx_k,$ for every $k,m\in\N$ and $x\in S$.}

{\bf Proof.} According to Theorem 3.3 we have \bea (V^m(x))_k &= &
(V^{m-1}(x))_k\left(1+\sum_{i=1}^{\x}a_{ki}(V^{m-1}(x))_i\right)\nonumber\\
&\leq & 2(V^{m-1}(x))_k\leq\cdots\leq 2^mx_k\nonumber \eea this is
the required relation.

{\bf Theorem 6.4.} {\it For every $m\in \N$ the sequence
$\{\V^m_n\}$ converges.}

{\bf Proof.} To show the convergence it is enough to prove that
$\{\V^m_n(x)\}$ is a Cauchy sequence for every $x\in S$. Without
loss of generality we may assume that $x_k>0$ for all $k\in\N$.

Let $\e>0$ be an arbitrary number. Since $x\in S$ there is  a
number $n_0\in\N$ such that (\ref{n0}) holds. Consider several
cases.

{\tt Case (i).} Let $1\leq k\leq n$ and $p\in\N$ be an arbitrary
number. For the sake of brevity we will denote \be\label{ab}
a^{(s)}_k=(V^s_{n]}(x))_k, \ \ \ b^{(s)}_k=(V^s_{n+p]}(x))_k, \ee
where $s,k\in\N$. Then from (\ref{comp2}),(\ref{vol4}) and
(\ref{ab}) we have \bea\label{fun}
|(\V^m_n(x))_k-(\V^m_{n+p}(x))_k|& =
&|a^{(m)}_k-b^{(m)}_k|\nonumber \\
&= &
\bigg|a^{(m-1)}_k\bigg(1+\sum_{i=1}^{\x}a^{n]}_{ki}a^{(m-1)}_i\bigg)
- b^{(m-1)}_k\bigg(1+\sum_{i=1}^{\x}a^{n+p]}_{ki}b^{(m-1)}_i\bigg)\bigg|\nonumber\\
&\leq &
|a^{(m-1)}_k-b^{(m-1)}_k|+\left|\sum_{i=1}^na^{n]}_{ki}\bigg(a^{(m-1)}_ka^{(m-1)}_i-b^{(m-1)}_kb^{(m-1)}_i\bigg)\right|\nonumber\\
& & + b^{(m-1)}_k\sum_{j=n+1}^{n+p}|a^{n+p]}_{kj}|b^{(m-1)}_j\nonumber \\
&\leq &
|a^{(m-1)}_k-b^{(m-1)}_k|+|a^{(m-1)}_k-b^{(m-1)}_k|\sum_{i=1}^n|a^{n]}_{ki}|a^{(m-1)}_i\nonumber
\\
& & +
b^{(m-1)}_k\sum_{i=1}^n|a^{n]}_{ki}|a^{(m-1)}_i-b^{(m-1)}_i|+
b^{(m-1)}_k\sum_{j=n+1}^{n+p}b^{(m-1)}_j\nonumber \\
&\leq & 2|a^{(m-1)}_k-b^{(m-1)}_k|+
b^{(m-1)}_k\sum_{i=1}^n|a^{(m-1)}_i-b^{(m-1)}_i|\nonumber\\
& & +b^{(m-1)}_k\sum_{j=n+1}^{n+p}b^{(m-1)}_j. \eea

Now we need the following

{\bf Lemma 6.5.} {\it For every $m\in\N$ the following inequality
holds \be\label{ab1}|a^{(m)}_k-b^{(m)}_k|\leq \a_m
x_k\sum_{j=n+1}^{n+p}x_j, \ee where
$$
\a_1=1, \ \a_m=\a_{m-1}(2+2^{m-1})+2^{2(m-1)}, \ \ \  \ m\geq 2.
$$}

{\bf Proof.} Let us firstly consider the case $m=1$. We have \bea
|a^{(1)}_k-b^{(1)}_k|& = &
|(V_{n]}(x))_k-(V_{n+p]}(x))_k|\nonumber \\
& = &
\left|x_k\bigg(\sum_{j=n+1}^{n+p}a^{n+p]}_{kj}x_j\bigg)\right|\leq
x_k\sum_{j=n+1}^{n+p}x_j.\nonumber \eea

This shows that $\a_1=1$.  Now assume that (\ref{ab1}) is valid
for $m-1$. Show that it is true for $m$. Indeed, it follows from
(\ref{fun}) and Lemma 6.3 that \bea|a^{(m)}_k-b^{(m)}_k|&\leq &
2\a_{m-1}x_k\sum_{j=n+1}^{n+p}x_j+\a_{m-1}2^{m-1}x_k\sum_{i=1}^nx_i\sum_{j=n+1}^{n+p}x_j\nonumber
\\
& & + 2^{2(m-1)}x_k\sum_{j=n+1}^{n+p}x_j\nonumber \\
&\leq &
(\a_{m-1}(2+2^{m-1})+2^{2(m-1)})x_k\sum_{j=n+1}^{n+p}x_j\nonumber
\eea which proves the lemma.

Now continue the proof of Theorem 6.4. According to Lemma 6.5 we
find that $ |(\V^m_n(x))_k-(\V^m_{n+p}(x))_k|<\e,$ for every
$n\geq n_0$.

{\tt Case (ii).} Let $n+1\leq k\leq n+p$. We have \bea\label{fun1}
|(\V^m_n(x))_k-(\V^m_{n+p}(x))_k|& = & \left|x_k-
b^{(m-1)}_k\bigg(1+\sum_{i=1}^{\x}a^{n+p]}_{ki}b^{(m-1)}_i\bigg)\right|\nonumber\\
&\leq &
|x_k-b^{(m-1)}_k|+b^{(m-1)}_k\sum_{i=1}^{n+p}b^{(m-1)}_i\nonumber\\
&\leq & |x_k-b^{(m-1)}_k|+2^{m-1}x_k.\eea

Now consider \bea |x_k-b^{(m-1)}_k|&\leq &
|x_k-b^{(m-2)}_k|+b^{(m-2)}_k\sum_{i=1}^{n+p}b^{(m-2)}_i\nonumber\\
&\leq&\cdots\leq
|x_k-b^{(1)}_k|+\sum_{j=1}^{m-2}b^{(j)}_k\sum_{i=1}^{n+p}b^{(j)}_i\nonumber\\
&\leq & x_k\sum_{i=1}^{n+p}x_i+x_k\sum_{j=1}^{m-2}2^{j}\leq
x_k\sum_{j=0}^{m-2}2^{j}.\nonumber \eea

Hence, from (\ref{fun1}) we infer that
$$
|(\V^m_n(x))_k-(\V^m_{n+p}(x))_k|\leq
\left(\sum_{j=0}^{m-1}2^{j}\right)\e
$$
for every $n\geq n_0$.

Now let $k\geq n+p$ then $(\V^m_n(x))_k=(\V^m_{n+p}(x))_k$.

Thus we have proved that  $\{\V^m_n(x)\}$ is a Cauchy sequence.
The limit of this sequence we  denote as $\W_m(x)$. The theorem is
proved.

From this theorem naturally arises a question: whether does  the
equality $\W_m=\V^m$ hold?

Before answer to this question we should  prove the following an
auxiliary

{\bf Lemma 6.6.} {\it Let $V$ be an arbitrary Volterra operator.
Then the following inequality holds
$$
\|V(x)-V(y)\|_1\leq 3\|x-y\|_1
$$
for every $x,y\in S$.}

{\bf Proof.} We have \bea
\|V(x)-V(y)\|_1&=&\sum_{k=1}^\x|(V(x))_k-(V(y))_k|\nonumber\\
&\leq &
\sum_{k=1}^\x\left(\bigg(1+\sum_{i=1}^\x|a_{ki}|x_i\bigg)|x_k-y_k|+x_k\sum_{i=1}^\x|a_{ki}||x_i-y_i|\right)\nonumber\\
&\leq &
2\sum_{k=1}^\x|x_k-y_k|+\sum_{i=1}^\x|x_i-y_i|=3\|x-y\|_1\nonumber
\eea

Lemma is proved.

{\bf Theorem 6.7.} {\it For every $m\in\N$ the equality
$\W_m=\V^m$ is valid.}

{\bf Proof.} Let $x\in S$ be an arbitrary element. Then given
$\e>0$ there is a number $n\in \N$ and $y\in S^{K_n}$ such that
$\|x-y\|_1<\e$. According to the compatibility condition
(\ref{comp1}) we have $\V(y)\in S^{K_n}$ and hence
$\V^m(y)=V_{n]}^m(y)$ therefore $\W_m(y)=\V^m(y)$. Using this we
have \be\label{conv0}
|(\W_m(x))_k-(\V^m(x))_k|\leq|(\W_m(x))_k-(\W_m(y))_k|+||(\V^m(x))_k-(\V^m(x))_k|
\ee for every $k\in\N$.

According to Theorem 6.4 we know that there is $n_0\in\N$ such
that \be\label{conv1}|(\W_m(x))_k-(\V^m_n(x))_k|<\e \ee for every
$n\geq n_0$.

Using Lemma 6.6 one gets \be\label{conv2}
|(\V^m_n(x))_k-(\V^m_n(y))_k|\leq\|\V^m_n(x)-\V^m_n(y)\|_1\leq
3^m\|x-y\|_1<3^m\e. \ee

It follows from  (\ref{conv1}),(\ref{conv2}) that
\bea\label{conv3} |(\W_m(x))_k-(\W_m(y))_k|&\leq
&|(\W_m(x))_k-(\V^m_n(x))_k|+|(\V^m_n(x))_k-(\V^m_n(y))_k|\nonumber\\
&&+|(\V^m_n(y))_k-(\W_m(y))_k|\leq (1+3^m)\e, \eea here we have
used the equality $\W_m(y)=\V^m_n(y)$.

Now again using Lemma 6.6 we find \be\label{conv4}
|(\V^m(x))_k-(\V^m(y))_k|\leq 3^m\|x-y\|_1<3^m\e. \ee

Consequently, the inequalities (\ref{conv3}),(\ref{conv4}) with
(\ref{conv0}) imply that
$$
|(\W_m(x))_k-(\V^m(x))_k|<(1+2\cdot 3^m)\e.
$$
As $\e$ has been an arbitrary, so this completes the proof.

This Theorem gives us some how to investigate limit behaviors of
infinite dimensional Volterra operators by means  of finite
dimensional ones. This would be a theme of our next
investigations.

{\bf Acknowledgements.}  The first named author (F.M.) thanks
NATO-TUBITAK for providing financial support and Harran University
for kind hospitality and providing all facilities. The work is
also partially supported by Grant $\Phi$-1.1.2 of Rep. Uzb.

\end{document}